\documentclass[12pt,a4paper]{amsart}

\usepackage{amsmath, amsfonts, xifthen, latexsym, amssymb, amsthm, amscd}

\usepackage[utf8]{inputenc}
\usepackage{graphicx}
\usepackage{url}

\usepackage{hyperref}
\hypersetup{colorlinks=true,citecolor=blue,filecolor=blue,linkcolor=blue,urlcolor=blue}
\usepackage[margin=1.4in]{geometry}

\newtheorem{theorem}{Theorem}
\newtheorem{lemma}{Lemma}[section]

\newtheorem{proposition}{Proposition}[section]

\newtheorem{corollary}{Corollary}[section]

\theoremstyle{definition}

\renewcommand{\phi}{{\varphi}}

\newcommand{\Free}{\operatorname{Free}}

\newcommand{\cC}{\mathcal{C}}

\newcommand\N{\mathbb N}

\newcommand{\actson}{\curvearrowright}

\newcommand{\cal}[1]{{\mathcal #1}}

\usepackage{etoolbox}
\newtoggle{no_cases}\toggletrue{no_cases}
\newcommand{\case}[2][]{\iftoggle{no_cases}{\left\{\begin{array}{ll}#2 & #1}{\\#2 & #1}\togglefalse{no_cases}}
\newcommand{\esac}{\end{array}\right.\toggletrue{no_cases}}

\newcommand{\Stab}{\operatorname{Stab}}\newcommand{\Sub}{\operatorname{Sub}}

\begin{document}
\title{On universal continuous actions on the Cantor set}
\author{Gábor Elek}


\thanks{The author was partially supported
by the ERC Consolidator Grant "Asymptotic invariants of discrete groups,
sparse graphs and locally symmetric spaces" No. 648017. }

\begin{abstract} Using the notion of proper Cantor colorings we prove the following theorem.
For any countably infinite group $\Gamma$, there exists a free continuous action $\zeta: \Gamma \actson \cC$ 
on the Cantor set, which is universal in the following sense:
for any free Borel action  $\alpha: \Gamma \actson X$ on the standard Borel space, there exists an injective Borel map
$\Theta_\alpha: X\to \cC$ such that $\Theta_\alpha\circ \alpha=\zeta \circ \Theta_\alpha$.
We  extend our theorem for (nonfree) Borel $(\Gamma,Z)$-actions, where $Z$ is a 
uniformly recurrent subgroup. 

\end{abstract}\maketitle

\setcounter{tocdepth}{2}
\section{Introduction}
In \cite{HM} Hjorth and Molberg proved that for any countable group $\Gamma$
there exists a free and continuous action of $\Gamma$ on the Cantor set $\cC$ that admits
an invariant probability measure. The first goal of this note is to show that one can say
much more.
\begin{theorem} \label{fotetel}
For any countably infinite group $\Gamma$, there exists a free continuous action $\zeta: \Gamma \actson \cC$ 
on the Cantor set such that
for any free Borel action  $\alpha: \Gamma \actson X$ on the standard Borel space, there exists an injective Borel map
$\Theta_\alpha: X\to \cC$, satisfying $\Theta_\alpha\circ \alpha=\zeta \circ \Theta_\alpha$.
\end{theorem}
\noindent
Therefore, if the free Borel action $\alpha:\Gamma\to X$ admits an invariant
probability  measure 
$\mu$  and
$\mu_{\Theta_\alpha}$ is the pushdown of $\mu$ onto $\cC$, then the probability measure $\mu_{\Theta_\alpha}$ is
invariant with respect to $\zeta$.
The theorem above will easily follow from the proposition below.
\begin{proposition} \label{free}
For any countably infinite group $\Gamma$ and for any free Borel action $\alpha:\Gamma\to X$ on the standard Borel space, there exists
an injective equivariant Borel map $\Theta'_\alpha: X \to \Free(\cC^\Gamma)$ ( where $\Free(\cC^\Gamma)$ is the free part of the
Bernoulli $\cal{C}$-shift space), such that the closure of the set $\Theta'_\alpha(X)$
is still in  $\Free(\cC^\Gamma)$.
\end{proposition}
\noindent
Note that Seward and Tucker-Drob \cite{ST} (see also \cite{Bern}) 
proved the following result:
For any free Borel action $\alpha:\Gamma\to X$ there exists a (not necessarily
injective)  equivariant Borel map $\Theta: X\to \Free(\{0,1\}^\Gamma)$
such that the closure of the set $\Theta(X)$ is still in $\Free(\{0,1\}^\Gamma)$. 
It seems however that the technique of \cite{ST} cannot easily be extended to
nonfree actions.

\noindent
So, let $\beta:\Gamma \to M$ be a continuous minimal nonfree action of a countable group 
on a compact metric space. Then we have a natural equivariant 
map $\Stab_\beta: M\to \Sub(\Gamma)$ (see Section \ref{uniform} for further details)
from our space $M$ to the compact space of subgroups of $\Gamma$, mapping each point $x\in M$ to its stabilizer subgroup.
Glasner and Weiss (Proposition 1.2 \cite{GW}) proved that the set $M_0$ of the 
continuity points
of the map $\Stab_\beta$ is a dense, invariant $G_\delta$ subset of $M$ and the closure
of $\Stab_\beta(M_0)$ in $\Sub(\Gamma)$ is a minimal closed invariant subset of
$\Sub(\Gamma)$, that is, a {\bf uniformly recurrent subgroup}. They asked if
for any uniformly recurrent subgroup $Z\in \Sub(\Gamma)$ there exists a
minimal continuous action $\beta:\Gamma \to M$ such that $M_0=M$. This question has been
answered in \cite{Elek} and \cite{MT}.
In Section \ref{uniform}
we will prove that one can answer the question of Glasner and Weiss in a uniform way. 
\begin{theorem} \label{masodik} Let $\Gamma$ be a countable group and $Z\in\Sub(\Gamma)$ be 
uniformly recurrent subgroup. Then there exists a continuous action $\zeta_Z:\Gamma\actson \cC$ such
that
\begin{itemize}
\item $\Stab_{\zeta_Z}: \cC\to \Sub(\Gamma)$ is continuous everywhere and $\Stab_{\zeta_Z}(\cC)=Z$.
\item For any Borel action $\alpha: \Gamma \actson X$  such that for any $x\in X$
the group $\Stab_\alpha(x)$ is in $Z$ (we call these actions $(\Gamma,Z)$-actions),
there exists an injective Borel map $\Theta_\alpha: X \to \cC$ such that $\Theta_\alpha\circ \alpha=\zeta_Z\circ \Theta_\alpha.$
\end{itemize}
\end{theorem}
\noindent
It was proved in Section 5 \cite{Elek} that there exist countable groups $\Gamma$ and uniformly recurrent subgroups
$Z\subset \Sub(\Gamma)$  such that no Borel $(\Gamma,Z)$-action admits an invariant probability measure. However,
we have the following nonfree analogue of the aforementioned result of Hjorth and Molberg.
\begin{corollary}
Let $\Gamma$ be a countable group and let $Z\subset \Sub(\Gamma)$ be a uniformly recurrent subgroup.
If there exists a Borel $(\Gamma,Z)$-action $\alpha:\Gamma\actson X$ that admits an invariant probability measure, then
there exists a continuous $(\Gamma,Z)$-action $\beta:\Gamma\actson \cC$ on the Cantor set admitting an invariant
probability measure such that
\begin{itemize}
\item The map $\Stab_\beta:\cC\to \Sub(\Gamma)$ is continuous everywhere and
\item $\Stab_\beta(\cC)=Z$.
\end{itemize}
\end{corollary}
\section{The proof of Theorem \ref{fotetel}}
Let $\Gamma$ be a countable infinite group and $\{\sigma_i\}^\infty_{i=1}$ be a generating system of $\Gamma$. Also for $n\geq 1$, let $\Gamma_n$ be the
subgroup of $\Gamma$ generated by the elements $\{\sigma_i\}^n_{i=1}$.
Let $\alpha:\Gamma\actson X$ be a Borel action of $\Gamma$ on the standard Borel space $X$.
We define a sequence $\{G_n\}^\infty_{n=1}$ of Borel graph structures on $X$ in the following way.
If $p,q\in X$, $p\neq q$, then let $(p,q)\in E(G_n)$ if there exists $1\leq i \leq n$ such that $\alpha(\sigma_i)(p)=q$ or $\alpha(\sigma_i)(q)=p$.
A Borel $\cal{C}$-coloring of $X$ is a Borel map $\phi:X\to \cC$, where $\cC=\{0,1\}^\N$ is the Cantor set. We say that $\phi$ is a {\bf proper} $\cC$-coloring
with respect to $\alpha:\Gamma\actson X$ if for any $r>0$ there exists $S_r>0$ such that
for any $p,q\in X$ $(\phi(p))_{S_r}\neq (\phi(q))_{S_r}$, provided that $0<d_{G_r}(p,q)\leq r$, where
\begin{itemize}
\item $d_{G_r}$ is the shortest path metric on the components of the Borel graph $G_r$.
\item For $\kappa\in \cC$ and $s>0$, $(\kappa)_s\in \{0,1\}^{[s]}$, denotes the projection of $\kappa$ onto its first $s$ coordinates.
\end{itemize}
\noindent
The Borel coloring $\phi$ is called separating if for any $p\neq q\in X$, $\phi(p)\neq \phi(q)$.
\begin{lemma} \label{proper} For any Borel action $\alpha:\Gamma\actson X$ there exists a separating, proper $\cC$-coloring with respect to $\alpha$.
\end{lemma}
\proof
First, for any $r\geq 1$ we construct
a new Borel graph $H_r$ of bounded vertex degree on $X$ such that $(p,q)\in E(H_r)$ if $0<d_{G_r}(p,q)\leq r$.
By the classical result of Kechris, Solecki and Todorcevic \cite{KST}, there exists an integer $m_r>0$ and a Borel coloring
$\psi_r:X\to \{0,1\}^{[m_r]}$ such that
$\psi_r(p)\neq \psi_r(q)$, whenever $p$ and $q$ are adjacent vertices in the Borel graph $H_r$.
Then $\phi_1(p)=\{\psi_1(p)\psi_2(p)\dots\}\in\{0,1\}^\N$, defines
a proper $\cC$-coloring of $X$ with respect to $\alpha$. Now we use the usual trick to obtain a separating coloring.
Let $\phi_2:X\to \cC$ be a Borel isomorphism. Then if $\phi_1(p)=\{a_1 a_2 a_3\dots \}$ and
$\phi_2(p)=\{b_1 b_2 b_3\dots \}$ let
$\phi(p)=\{a_1b_1a_2b_2\dots\}$. Clearly, $\phi$ is a separating, proper $\cC$-coloring with respect to $\alpha$. \qed
\vskip 0.1in
\noindent
Now we prove Proposition \ref{free}. Let $\alpha:\Gamma \actson X$ be a free Borel action and let $\phi:X\to \cC$ be a separating, proper $\cC$-coloring.
Consider the Bernoulli shift $\cC^\Gamma$ with the natural
left action
$$L_\delta(\rho)(\gamma)=\rho(\gamma\delta)\,.$$
\noindent
The map $\Theta'_\alpha:X\to \cC^\Gamma$ is defined as usual by
$$\Theta'_\alpha(x)(\gamma)=\phi(\alpha(\gamma)(x))\,.$$
\noindent
Clearly, $\Theta'_\alpha$ is Borel and $\Gamma$-equivariant and since
$\phi$ is separating $\Theta'_\alpha$ is injective as well. Also, $\Theta'_\alpha(X)\subset \Free(\cC^\Gamma)$.
We need to show that if
$\rho:\Gamma\to\cC$ is in the closure of $\Theta'_\alpha(X)$, then
$L_\gamma(\rho)\neq \rho$, whenever $\gamma\neq 1_\Gamma$.
It is enough to see that if $\gamma\neq e_\Gamma$, then
$\rho(\gamma)\neq \rho(e_\Gamma)$.
Let $\lim_{n\to\infty} \Theta'_\alpha(x_n)=\rho\in\cC^\Gamma$.
Also, let $r>0$ such that
\begin{itemize}
\item $\gamma\in\Gamma_r.$
\item $d_{\Gamma_r}(e_\Gamma,\gamma)\leq r$, where $d_{\Gamma_r}$ is the shortest path metric on the right Cayley graph $\mbox{Cay}(\Gamma_r,\{\sigma_i\}^r_{i=1})$.
\end{itemize}
\noindent
Thus, for any $n\geq 1$,
$$d_{G_r}(\alpha(\gamma) (x_n), x_n)\leq r\,.$$
\noindent
Hence, for any $n\geq 1$,
$$\left(\phi(\alpha(\gamma)(x_n))\right)_r\neq (\phi(x_n))_r\,.$$
\noindent
Since $\rho(\gamma)=\lim_{n\to\infty} \phi(\alpha(\gamma)(x_n))$ and
$\rho(e_\Gamma)=\lim_{n\to\infty} \phi(x_n)$, we have that
$\rho(\gamma)\neq \rho(e_\Gamma)\,.$ Hence our proposition follows. \qed
\vskip 0.1in
\noindent
Now we prove Theorem \ref{fotetel}.
Let $\iota:\Gamma\actson \Free(\cC^\Gamma)$ be
the natural free action and $\Theta'_\iota:\Free(\cC^\Gamma)\to \Free(\cC^\Gamma)$ be the
injective $\Gamma$-equivariant Borel map defined in Proposition \ref{free}.
Let $W$ be the closure of $\Theta'_\iota(\Free(\cC^\Gamma))$ in $\cC^\Gamma$. Thus, $W\subset \Free(\cC^\Gamma)$ is a closed invariant subset.
Then, $W$ can be written as the disjoint union of $W_1$ and $W_2$, where $W_1$ is a countable invariant subset and $W_2$ is a closed invariant subset
homeomorphic to the Cantor set.
Let $\alpha:\Gamma\actson X$ be a free Borel action and $\Theta_\alpha':X\to \Free(\cC^\Gamma)$ be the 
injective $\Gamma$-equivariant  Borel map defined in  Proposition \ref{free}.
Let $\tilde{\Theta}_\alpha=\Theta'_\iota \circ \Theta'_\alpha.$
Clearly, $\tilde{\Theta}_\alpha:X\to W$ is an injective $\Gamma$-equivariant Borel map. Now we define the injective $\Gamma$-equivariant Borel map
$\Theta_\alpha:X\to W_2$ by modifying $\tilde{\Theta}_\alpha$ on countably many orbits. Since $W_2$ is homeomorphic to $\cC$ our Theorem follows. \qed

\section{Uniformly recurrent subgroups} \label{uniform}
Let $\Gamma$ be a countable group and let $\Sub(\Gamma)$ be the  space of subgroups of $\Gamma$ \cite{GW}. Then $\Sub(\Gamma)$ is
a compact, metrizable space and  conjugations define a continuous action $c:\Gamma\actson \Sub(\Gamma)$. Let $Z\subset \Sub(\Gamma)$ be
a uniformly recurrent subgroup as in the Introduction.
We define the Bernoulli shift space $\cC^Z$ of $Z$ in the following way.
Let $$\cC^Z=\cup_{H\in Z} \cal{F}(H)\,,$$
where $\cal{F}(H)$ is the set of maps $\rho:\Gamma/H\to \cC$ from the right coset space of $H$ to the Cantor set.
\noindent
The action of $\Gamma$ on $\cC^Z$ is defined as follows.
\begin{itemize}
\item If $\rho\in \cal{F}(H)$ then
$L_\delta(\rho)\in \cal{F}(\delta H \delta^{-1})$ and
\item $L_\delta(\rho)(\delta H \delta^{-1}\gamma)=\rho(H\gamma\delta)\,.$
\end{itemize}
\begin{lemma}
$L:\Gamma\to \mbox{Homeo}(\cC^Z)$ is a homomorphism.
\end{lemma}
\proof
We need to show that if
$\rho:\Gamma/H\to \cC$ and $\delta_1,\delta_2\in \Gamma$, then
$$L_{\delta_1 }(L_{\delta_2}(\rho))=L_{\delta_1\delta_2}(\rho)\,.$$
\noindent
Observe that $$L_{\delta_1} (L_{\delta_2}(\rho))\in\cal{F}(\delta_1\delta_2H\delta_2^{-1}\delta^{-1}_1)\,\,\mbox{and}\,\,
L_{\delta_1\delta_2}(\rho)\in \cal{F}(\delta_1\delta_2H\delta_2^{-1}\delta^{-1}_1).$$
Now
$$L_{\delta_1}(L_{\delta_2}(\rho))(\delta_1\delta_2H\delta_2^{-1}\delta^{-1}_1\gamma)= $$ $$
=L_{\delta_2}(\rho)(\delta_2 H \delta^{-1}_2\gamma\delta_1)=\rho (H\gamma\delta_1\delta_2)=
L_{\delta_1\delta_2}(\rho)(\delta_1\delta_2H\delta_2^{-1}\delta^{-1}_1\gamma)\,,$$
\noindent
hence our lemma follows. \qed
\vskip 0.1in
\noindent
We can equip $\cC^Z$ with a compact metric structure $d$ such that
$(\cC^Z,d)$ is homeomorphic to the Cantor set and the $\Gamma$-action above is continuous. 
Let $\rho_1:\Gamma/H_1\to \cC$, $\rho_2:\Gamma/H_2\to \cC$ be elements of $\cC^Z$. We say that
$\rho_1$ and $\rho_2$ are $n$-equivalent, $\rho_1 \equiv_n \rho_2$ if
\begin{itemize}
\item For any $\gamma\in \Gamma$, $d_{\Gamma_n}(e_\Gamma,\gamma)\leq n$,
$\gamma\in H_1$ if and only if $\gamma\in H_2$. 
\item For any $\gamma\in \Gamma$, $d_{\Gamma_n}(e_\Gamma,\gamma)\leq n$,
$$(\rho_1(H_1\gamma))_n=(\rho_2(H_2\gamma))_n\,.$$
\end{itemize}
\noindent
Then we define $d(\rho_1,\rho_2):=\frac{1}{2^n}$ whenever
$\rho_1\equiv_n \rho_2$ and $\rho_1 \not\equiv_{n+1} \rho_2$.
Let $\{H_n\}^\infty_{n=1}, H\in \Sub(\Gamma)$, $\rho\in\cal{F}(H)$ and for
any $n\geq 1$ let $\rho_n\in\cal{F}(H_n)$.
Observe that $\{\rho_n\}^\infty_{n=1}\to \rho$ in the $d$-metric if and only if
\begin{itemize}
\item $H_n\to H$ in the compact space $\Sub(\Gamma)$ and
\item for any $\gamma\in\Gamma$, $\rho_n(H_n\gamma)\to \rho(H\gamma)\,.$
\end{itemize}
\noindent
We can define $\Free(\cC^Z)$ in the usual way.
We have that $\rho:\Gamma/H\to \cC\in \Free(\cC^Z)$ if
$L_\delta(\rho)\neq \rho$ for any $\delta\notin H$.
Clearly, if for any $\delta\notin H$, $\rho(H)\neq \rho(H\delta)$, then $\rho\in \Free(\cC^Z)$.
Now let $\alpha:\Gamma\actson X$ be a Borel $(\Gamma,Z)$-action and $\phi:X\to \cC$ be a separating, proper $\cC$-coloring with respect
to $\alpha$.
We define $\Theta^\phi_\alpha:X\to \cC^Z$ as follows.
\begin{itemize}
\item $\Theta^\phi_\alpha(x)\in \cal{F}(H)$, where $H=\Stab_\alpha(x)$.
\item $\Theta^\phi_\alpha(x)(H\gamma)=\phi(\alpha(\gamma)(x))\,.$
\end{itemize}
\noindent
Clearly, $\Theta^\phi_\alpha:X\to \cC^Z$ is an injective Borel map and $\Theta^\phi_\alpha(X)\subset \Free(\cC^Z)$.
\begin{lemma} 
The map $\Theta_\alpha^\phi:X\to \Free(\cC^Z)$ is $\Gamma$-equivariant.
\end{lemma}
\proof
Let $\delta\in\Gamma$. Then
$$L_\delta(\Theta^\phi_\alpha(x))(\delta H \delta^{-1}\gamma )=\Theta^\phi_\alpha(x)(H\gamma\delta)=
\phi(\alpha(\gamma\delta)(x))\,.$$
\noindent
On the other hand,
$$\Theta^\phi_\alpha(\alpha(\delta)(x))(\delta H\delta^{-1}\gamma)=\phi(\alpha(\delta)\alpha(\gamma)(x))=
\phi(\alpha(\gamma\delta))(x))\,.\quad\qed $$
\noindent
Now we prove the nonfree analogue of Proposition \ref{free}.
\begin{proposition} \label{free2}
For any countably infinite group $\Gamma$ and for any free Borel 
$(\Gamma,Z)$-action $\alpha:\Gamma\actson  X$, there exists
an injective equivariant Borel map $\Theta'_\alpha: X \to \Free(\cC^Z)$, such that the closure of the set $\Theta'_\alpha(X)$
is still in  $\Free(\cC^Z)$.
\end{proposition}
\proof Let $\phi$ and $\Theta^\phi_\alpha$ be as above.
Let $\{x_n\}^\infty_{n=1}\subset X$ such that
$$\lim_{n\to \infty}\Theta^\phi_\alpha(x_n)=\rho\in\cal{F}(H)$$ \noindent and $\delta \notin H$. We need to show that $\rho\in \Free(\cC^Z)$. 
Observe that \\ $\{\Stab_\alpha(x_n)\}^\infty_{n=1}\to H$ in $\Sub(\Gamma)$. 
Hence, there exists $N>0$ such that $\delta\notin H_n$ if $n\geq N$. By properness, there exists $m>0$ such
that for all $n\geq N$
$$(\phi(\alpha(\delta) (x_n)))_m\neq (\phi(x_n))_m\,.$$
\noindent
Since, $\lim_{n\to\infty} \phi(\alpha(\delta)( x_n))=\rho(H\delta)$ and $\lim_{n\to\infty} \phi(x_n)=\rho(H)$
we have that $\rho(H)\neq \rho(H\delta)$. Therefore $\rho\in \Free(\cC^Z)$. \qed
\vskip 0.1in
\noindent
Now Theorem \ref{masodik} follows from Proposition \ref{free2} exactly
the same way as Theorem \ref{fotetel} follows from Proposition \ref{free}. \qed

\end{document}